# EXISTENCE AND NUMBER OF SOLUTIONS OF DIOPHANTINE QUADRATIC EQUATIONS WITH TWO UNKNOWNS IN $\mathbb{Z}$ AND $\mathbb{N}$

Florentin Smarandache
University of New Mexico
200 College Road
Gallup, NM 87301, USA
E-mail: smarand@unm.edu

**Abstract:** In this short note we study the existence and number of solutions in the set of integers (Z) and in the set of natural numbers (N) of Diophantine equations of second degree with two unknowns of the general form $ax^2 - by^2 = c$.

**Property 1:** The equation $x^2 - y^2 = c$ admits integer solutions if and only if $c$ belongs to $4\mathbb{Z}$ or is odd.

*Proof:* The equation $(x-y)(x+y) = c$ admits solutions in $\mathbb{Z}$ iff there exist $c_1$ and $c_2$ in $\mathbb{Z}$ such that $x - y = c_1$, $x + y = c_2$, and $c_1 c_2 = c$. Therefore $x = \dfrac{c_1 + c_2}{2}$ and $y = \dfrac{c_2 - c_1}{2}$.

But $x$ and $y$ are integers if and only if $c_1 + c_2 \in 2\mathbb{Z}$, i.e.:

1) or $c_1$ and $c_2$ are odd, then $c$ is odd (and reciprocally).
2) or $c_1$ and $c_2$ are even, then $c \in 4\mathbb{Z}$.

Reciprocally, if $c \in 4\mathbb{Z}$, then we can decompose up $c$ into two even factors $c_1$ and $c_2$, such that $c_1 c_2 = c$.

**Remark 1:**

Property 1 is true also for solving in $\mathbb{N}$, because we can suppose $c \geq 0$ {in the contrary case, we can multiply the equation by (-1)}, and we can suppose $c_2 \geq c_1 \geq 0$, from which $x \geq 0$ and $y \geq 0$.

**Property 2:** The equation $x^2 - dy^2 = c^2$ (where $d$ is not a perfect square) admits an infinity of solutions in $\mathbb{N}$.

*Proof:* Let's consider $x = ck_1$, $k_1 \in \mathbb{N}$ and $y = ck_2$, $k_2 \in \mathbb{N}$, $c \in \mathbb{N}$. It results that $k_1^2 - dk_2^2 = 1$, which we can recognize as being the Pell-Fermat's equation, which admits an infinity of solutions in $\mathbb{N}$, $(u_n, v_n)$. Therefore $x_n = cu_n$, $y_n = cv_n$ constitute an infinity of natural solutions for our equation.

**Property 3**: The equation $ax^2 - by^2 = c$, $c \neq 0$, where $ab = k^2$, ($k \in \mathbb{Z}$), admits a finite number of natural solutions.

*Proof:* We can consider $a$, $b$, $c$ as positive numbers, otherwise, we can multiply the equation by (-1) and we can rename the variables.

Let us multiply the equation by $a$, then we will have:

$z^2 - t^2 = d$ with $z = ax \in \mathbb{N}$, $t = ky \in \mathbb{N}$ and $d = ac > 0$. (1)

We will solve it as in property 1, which gives $z$ and $t$.
But in (1) there is a finite number of natural solutions, because there is a finite number of integer divisors for a number in $\mathbb{N}^*$. Because the pairs $(z,t)$ are in a limited number, it results that the pairs $(z/a, t/k)$ also are in a limited number, and the same for the pairs $(x, y)$.

**Property 4:** If $ax^2 - by^2 = c$, where $ab \neq k^2$ ($k \in \mathbb{Z}$) admits a particular nontrivial solution in $\mathbb{N}$, then it admits an infinity of solutions in $\mathbb{N}$.

*Proof:* Let's consider:

$$\begin{cases} x_n = x_0 u_n + b y_0 v_n \\ y_n = y_0 u_n + a x_0 v_n \end{cases} \quad (n \in \mathbb{N}) \qquad (2)$$

where $(x_0, y_0)$ is the particular natural solution for the initial equation, and $(u_n, v_n)_{n \in \mathbb{N}}$ is the general natural solution for the equation $u^2 - abv^2 = 1$, called the solution Pell, which admits an infinity of solutions.

Then $ax_n^2 - by_n^2 = (ax_0^2 - by_0^2)(u_n^2 - abv_n^2) = c$.

Therefore (2) verifies the initial equation.

[1982]